\documentclass[preprint,12pt]{elsarticle}

\usepackage{amsmath,amssymb,mathdots}
\usepackage{bm,upgreek,nicefrac}
\usepackage[mathscr]{eucal}
\usepackage{xcolor}
\usepackage{tikz,accents,subcaption}

\newcommand{\bi}[1]{\emph{\textbf{#1}}}

\newcommand{\undb}[1]{\underaccent{\bar}{#1}}
\newcommand{\undt}[1]{\underaccent{\tilde}{#1}}
\newcommand{\BL}{\bar{\mathcal{L}}^{(1,2)}_A}

\newcommand{\UBL}{\undb{\mathcal{L}}^{(1,2)}_A}

\newtheorem{thm}{Theorem}[section]
\newtheorem{cor}[thm]{Corollary}
\newtheorem{lem}[thm]{Lemma}

\newtheorem{exmp}{Example}[section]
\newproof{pf}{Proof}

\journal{Syst. Control Lett. }

\begin{document}
\begin{frontmatter}

\title{Laplacian Controllability of Interconnected Graphs}

\author{Shun-Pin Hsu\corref{cor1}}
\cortext[cor1]{corresponding author}
\ead{shsu@nchu.edu.tw}


\address{Department of Electrical Engineering, National Chung Hsing University\\
250, Kuo-Kuang Rd., Taichung 402, Taiwan}

\begin{abstract}
In this work we consider the Laplacian controllability of a graph constructed by interconnecting a finite number of single-input Laplacian controllable graphs. We first study the interconnection realized by the composite graph of two connected simple graphs called the structure graph and the cell graph. Suppose the cell graph is Laplacian controllable by an input connected to some special vertex called the composite vertex. The composite graph is constructed by interconnecting all cell graphs through the composite vertices which alone form the structure graph. We then show that the structure graph is Laplacian controllable by an input connected to some vertex of the graph if and only if the composite graph is Laplacian controllable by that input connected to that composite vertex. In the second part of the paper, we view a path as a graph generated by interconnecting a finite number of two-vertex antiregular graph, and possibly connected to a one-vertex path, where the two vertices of the antiregular graph are interpreted as the terminal vertex (or the dominating vertex) and the degree-repeating vertex. We show that with a similar connecting scheme the single-input Laplacian controllability is preserved if we increase the number of vertices of the antiregular graph and that of the path. Numerical examples are presented to illustrate our results.
\end{abstract}

\end{frontmatter}

\section{Introduction}
Modern control, communication and networking technologies admit their integrations to innovative systems composed of the units such as robots or drones~\cite{mes10}. However, the common property that each unit communicates only with other reachable units causes system-wide problems such as the synchronization, coordination and formation~\cite{lin05,ji09}. Many of the issues have been partially addressed under the formulation of multiagent systems following the leader-follower dynamics. A more fundamental topic closely related to the effective operation of the networked system is on its controllability. This topic has been a main theme of many studies in system and control sciences since Tanner's conference paper was published more than a decade ago~\cite{tan04}. This pioneering paper investigated the interplay between control and communication in a networked system defined on a graph and formulated a nontrivial new problem that inspired many research works in theses areas. In the general setting, a linear and time-invariant system observing the consensus policy is proposed to simulate the information propagation across a networked system. Each state variable in the system changes with the difference between its value and that of other variables interacting with it. All state variables become steady as these variables reach a common value, or the consensus. In this scenario an autonomous system evolves according to the Laplacian dynamics is induced~\cite{agu15}. The controllability issues are raised when input signals are used to maneuver the linear system. Testing the controllability of the system is a classical problem. Several solutions including Kalman's rank test are available~\cite{chen99}. Nevertheless, a major challenge to the controllability issue in a networked system is the numerical imprecision caused by the large scale of the system. A standard approach is to adopt the controllability test proposed by Popov, Belevitch and Hautus. This test sheds some light on the relation between Laplacian eigenspaces and Laplacian controllability of the graph modeling the networked system. Since specific Laplacian eigenspaces can be inferred from the connection topology of the graph, partitioning schemes such as the equitable partition and the almost equitable partition~\cite{car07} were proposed to identify certain symmetry property that leads to Laplacian uncontrollability ~\cite{rah09,eg12}. These schemes helped to find a class of uncontrollable graphs, but failed to recognize a controllable one. This failure is partially compensated by a combined use of the distance partition and almost equitable partition schemes to bound the dimension of controllable subspace of state variables~\cite{zha14}. This expedient is simple and feasible for a large-scale system, but it offers a tight bound only in  the extreme case such as a path. In fact, how to use the minimum number of controllers to drive a general networked system to achieve an arbitrary state is a challenging task. The solved cases are restricted to the networks following specific connection structures, including the paths~\cite{par12}, multi-chain~\cite{hsu17a}, antiregular graphs~\cite{agu15b}, grids~\cite{not13}, circulant graphs~\cite{nab13}, and complete graphs~\cite{zhang11}. In these examples Laplacian spectra and eigenspaces of the graphs are fully available and thus  their controllability analyses are tractable.

In the previous study, we proposed a sufficient condition for a graph to be controllable by one input~\cite{hsu17b}. The breakthrough of the result is showing that a partial knowledge of the Laplacian spectrum and eigenspaces might be enough to identify a Laplacian controllable graph. The conclusion applies to a special graph family in which each graph is constructed by interconnecting a path and an antiregular graph. This work follows the line of our previous study and investigates a richer class of graphs constructed by interconnecting a finite number of single-input Laplacian controllable graphs. We first consider two connected simple graphs, called the structure graph and cell graph, where the cell graph is Laplacian controllable by an input connected to some vertex, called the composite vertex, of the graph. We then interconnect all cell graphs to form a composite graph through the composite vertices such that composite vertices alone form the structure graph. We show that the structure graph is Laplacian controllable by an input connected to some vertex of the graph, if and only if the composite graph is Laplacian controllable by that input connected to that (composite) vertex. Furthermore, we extend our result by showing that using different vertices to interconnect the cell graphs might also lead to the single-input Laplacian controllability of the resulting graph.

The contributions of our results are to propose an innovative and simple method to generate a huge class of single-input Laplacian controllable graphs whose controllability properties cannot be inferred from existing methods. Moreover, our results provide potential solutions for the design problem that asks to construct a $k$-vertex Laplacian controllable graph whose edge related parameters such as the diameter and maximum degree are subject to constraints. Specifically, our first result shows that the graph can be generated from the composite graph of a $k_1$-vertex graph and a $k_2$-vertex graph where $k_1k_2=k$ as long as each of these two graphs are Laplacian controllable by an input connected to some vertex of the graph. Our second result provides another solution to generate the graph by  interconnecting $c$ antiregular graphs, each with $k_a$ vertices, and one $k_p$-vertex path, where $ck_a+k_p=k$ and $c,k_a,k_p$ are nonnegative integers (Note that in~\cite{hsu17b} $c$ is restricted to $1$ and the interconnecting scheme is different). As $k$ is large, our results become very helpful in providing a large set of potential solutions for the graph designers to determine the network topology of the system that meets the requirements.

The rest of the paper is organized as follows. In the second section we recapitulate the essential graph-theoretical notations and concepts as well as related control theories. In the third section we present our main results on the schemes to interconnect graphs while preserving the Lpalacian controllability. The paper is concluded in Section~4
where possible extensions of our results are discussed.

\section{Preliminaries}
We start this section by defining or reviewing notations and concepts to be used throughout the paper. Let $\bi{e}_i$ be the zero vector except its $i$th entry being $1$. The floor and ceiling functions of $x$, written as $\lfloor x\rfloor$ and $\lceil x\rceil$, are the largest integer not greater than $x$ and the smallest integer not less than $x$, respectively. The identity matrix of order $k$ is $I_k$. Occasionally we drop the subscript $k$ for simplicity as the context is clear. The difference of  sets $S_1$ and $S_2$ is $S_1\setminus S_2$, meaning that $\{s|s\in S_1, s\notin S_2\}$. If $\lambda$ is an eigenvalue of $P$ and $\bi{v}$ is an eigenvector corresponding to $\lambda$, we say $(\lambda,\bi{v})$ is an eigenpair of $P$. If $A_1$ and $A_2$ are square matrices of the same order, $A_1$ is called a modal matrix of $A_2$ if the columns of $A_1$ are independent eigenvectors of $A_2$. The Kronecker product of matrices $P_1$ and $P_2$ are written as $P_1\otimes P_2$. If $V=\{1,2,\cdots,k\}$ and $E$ is a subset of $\{\,(v_1,v_2)\,|\,v_1,v_2\in V\}$, then $(V,E)$ describes a $k$-vertex graph, or a graph on $k$ vertices, where $V$ and $E$ are called the vertex set and edge set respectively, of the graph. In this work we consider only the connected simple graphs to avoid apparent uncontrollability and to simplify the interconnection pattern. The algebraic aspects of such graphs can be seen, for example, in~\cite{god01,biggs93}. We say $v_1$ and $v_2$ are neighbors if $v_1,v_2\in V$ and $(v_1,v_2)\in E$. The neighbor set of the vertex $v$ is $\mathcal{N}_v:=\{u\,|\,(v,u)\in E\}$. The degree, or valency, of vertex $v$ is defined as $|\mathcal{N}_v|$, the number of elements in $\mathcal{N}_v$. In a connected simple graph defined by $(V,E)$, at least two vertices share the same degree. The vertex $v$ is called a terminal vertex if $|\mathcal{N}_v|=1$, it is called a dominating vertex if $|\mathcal{N}_v|=|V|-1$, and is called a degree-repeating vertex if its degree is not unique among those of vertices in the graph. With appropriate vertex numbering, we can assign $d_i$ to be the degree of the $i$th vertex of a $k$-vertex graph such that $d_i\geq d_{i+1}$ for each $i\in\{1,2,\cdots,k-1\}$. The sequence $\bi{d}:=(\,d_1,d_2,\cdots,d_k\,)$ is called the degree sequence of the $k$-vertex graph. The trace of $\bi{d}$ is $\tau_{\bi{d}}:=|\,j:d_j\geq j\,|$. The conjugate of $\bi{d}$ is $\bi{d}^*:=\left(\,d_1^*,d_2^*,\cdots,d_k^*\,\right)$  where $d_i^*=|\,j:d_j\geq i\,|$. Given an arbitrary sequence $\bi{d}=(\,d_1,d_2,\cdots,d_k\, )$ of nonnegative integers, $\bi{d}$  is called graphical if there exists a $k$-vertex graph whose degree sequence is $\bi{d}$. It was shown~\cite[p.72]{mol12} that the necessary and sufficient condition for $\bi{d}$ to be graphical is that
\begin{equation}\label{lem:gra}
\sum_{i=1}^j(d_i+1)\leq\sum_{i=1}^jd_i^*, \quad\mbox{$\forall j\in\{1,2,\cdots,\tau_\bi{d}\}$}.
\end{equation}
Suppose $(V,E)$ determines a connected simple graph and its degree sequence is $\bi{d}=(\,d_1,d_2,\cdots,d_k\, )$. The Laplacian matrix $\mathcal{L}$ of the graph is defined as
\begin{equation*}
\mathcal{L}:=\mathcal{D}-\mathcal{A}
\end{equation*}
where $\mathcal{D}$ is a diagonal matrix whose $i$th diagonal term is $d_i$, and $\mathcal{A}$ is a binary matrix whose $(i,j)$th element is $1$ if $(i,j)\in E$ and is $0$ otherwise. The eigenvalues and eigenvectors of $\mathcal{L}$
are called the Laplacian eigenvalues and Laplacian eigenvectors, respectively, of the graph. An apparent eigenpair of $\mathcal{L}$ is $(0,\textbf{1})$ where $\textbf{1}$ is a vector of $1$'s. Many properties concerning the eigenvalues and eigenvectors of $\mathcal{L}$ were summarized in~\cite{mer98}. An interesting result was recently proposed that the spectrum of $\mathcal{L}$ is majorized by the conjugate $\bi{d}^*$ of the degree sequence of the graph, namely,
\begin{equation}\label{lem:GroMer}
\sum_{i=1}^t\ell_{k-i+1}\le\sum_{i=1}^td_i^*\quad\forall t\in\{1,2,\cdots,k\}
\end{equation}
where $\ell_i$ is the $i$th smallest Laplacian eigenvalue of the graph. If the equality in~(\ref{lem:gra}) holds, then $\bi{d}$ determines uniquely a \emph{threshold graph} or \emph{maximal graph}. It was proved  in~\cite{mer94} that the equality in~(\ref{lem:gra}) implies the equality in~(\ref{lem:GroMer}). Thus the eigenvalues of a threshold graph is readily available from its degree sequence. In fact, threshold graphs admit many different definitions from the version used here, e.g., the definition based on the Ferrers-Sylvester diagram~\cite[p.70]{mol12}, and the one on constructing the graph using join and union operations only~\cite{bap13}. An antiregular graph is a connected simple graph that has exactly one pair of vertices sharing the same degree~\cite{mer03}. An antiregular graph turns out to be a special threshold graph and thus enjoys many excellent Laplacian eigenpair properties.

A linear and time-invariant system that has $k$ state variables and evolves according to the consensus policy
has the following form:
\begin{equation}\label{eq:cons}
\dot{x}_i=-\sum_{j\in\mathcal{N}_i}(x_i-x_j),
\end{equation}
where $\mathcal{N}_i$ is a subset of the set of state variables. This system can be readily used to model a dynamic networked system in which each vertex state interacts with its neighboring vertex states and expects to reach a common state as time goes by. As the model is used in this context, the evolution of the linear system follows the Laplacian dynamics~\cite[p.1613]{agu15} written as
\begin{equation}\label{eq:automcons}
\dot{\bi{x}}=-\mathcal{L}\bi{x}
\end{equation}
where $\mathcal{L}$ is the Laplacian matrix of the graph describing the networked system. To control the autonomous system in~(\ref{eq:automcons}), we can apply $p$ control inputs $\bi{u}(t)=[\,u_1(t)\,u_2(t)\,\cdots u_p(t)\,]^T$ via a binary control coefficient matrix $B\in\{0,1\}^{k\times p}$ so that
\begin{equation}\label{eq:mcons}
\dot{\bi{x}}=-\mathcal{L}\bi{x}+B\bi{u}(t)
\end{equation}
where the $(i,j)$th element of $B$ is
$1$ if vertex $i$ is connected to input $u_j(t)$, and is $0$ otherwise. For simplicity, we use the notation $(\mathcal{L},B)$ to represent the controlled graph model of the dynamic system in~(\ref{eq:mcons}). In the single-input system $p$ is $1$ and the specific notation $(\mathcal{L},\bi{b})$ is used. We say a graph is single-input Laplacian controllable if its corresponding $(\mathcal{L},\bi{b})$ is controllable. The focus of this paper is on the controllability of $(\mathcal{L},\bi{b})$ where $\mathcal{L}$ is the Laplacian matrix of a graph constructed by interconnecting a finite sequence of single-input Laplacian controllable graphs. Apparently, how the sequence of graphs are interconnected and how the input is applied affect the Laplacian controllability of the resulting graph. We first consider a scheme called the composite graph of two graphs and present a necessary and sufficient condition for the composite graph to be single-input Laplacian controllable. Following this result we show for a special case that the scheme of a composite graph that relies on specific vertices for interconnection can be modified to allow more flexibility while preserving the single-input Laplacian controllability. Before presenting our main results, we mention a version of
the classical Popov-Belevitch-Hautus (PBH) test on which our studies are based.
 \begin{thm}\label{thm:ns}~\cite[p.145]{chen99} A graph is Laplacian controllable if and only if the graph does not have a Laplacian eigenvector orthogonal to the column space of the control coefficient matrix.
 \end{thm}

\section{Main Results}
\subsection{the composite graph}
Suppose $\mathcal{L}^{(1)}$ and $\mathcal{L}^{(2)}$ are the Laplacian matrices of connected simple graphs $\mathbb{G}_1$ and $\mathbb{G}_2$, with $k_1$ and $k_2$ vertices respectively. If we interconnect $k_1$ copies of $\mathbb{G}_2$ via the $s$th vertex of each $\mathbb{G}_2$ such that these interconnecting vertices alone form $\mathbb{G}_1$, the resulting graph, written as
$\mathcal{G}_1(\mathbb{G}_2,s)$, is called the composite graph of $\mathbb{G}_1$ by $\mathbb{G}_2$ via vertex $s$. The graphs $\mathbb{G}_1$ and $\mathbb{G}_2$ are called the structure graph and cell graph respectively of $\mathcal{G}_1(\mathbb{G}_2,s)$, and the vertex $s$ in $\mathbb{G}_2$ the composite vertex. Let the Laplacian matrix of $\mathcal{G}_1(\mathbb{G}_2,s)$ be $\mathcal{L}_s^{(1,2)}$. We thus have
\begin{equation}
\mathcal{L}^{(1,2)}_s=I_{k_1}\otimes\mathcal{L}^{(2)}+\mathcal{L}^{(1)}\otimes\bi{e}_s\bi{e}_s^T.
\end{equation}
In the sequel we present several matrix properties that are closely related to the Laplacian eigenspace of the composite graph.

\begin{lem}\label{lem:indc} If $A$ is a square matrix of order $k$ and the $i$th entry of every eigenvector of $A$ is nonzero for some $i\in\{1,2,\cdots,k\}$, then all but the $i$th column of $A-\alpha I$ are independent for any $\alpha$.
\end{lem}
\begin{pf}
If for some $i\in\{1,2,\cdots,k\}$ the $i$th entry of every eigenvector of $A$ is nonzero, then every eigenvalue of $A$ is not repeated. In case $\alpha$ is not an eigenvalue of $A$, $A-\alpha I$ is nonsingular and thus all columns of $A-\alpha I$ are independent. If $\alpha$ is an eigenvalue of $A$, the nonzero entry in the eigenvector implies that the $i$th column of $A-\alpha I$ is a linear combination of the remaining columns. If all but the $i$th column of $A-\alpha I$ are dependent, then the rank of $A-\alpha I$ is at most $k-2$ and $\alpha$ is an eigenvalue with its algebraic multiplicity at least $2$, a contradiction.
\end{pf}

\begin{lem}\label{lem:nsing}
Let $\lambda$ be an eigenvalue of $\mathcal{G}_1(\mathbb{G}_2,s_2)$. If the $s_i$th entry of every eigenvector of $\mathcal{L}^{(i)}$ is nonzero for every $i$ in $\{1,2\}$, the matrix obtained by removing the $s_2$th row and $s_2$th column of $\mathcal{L}^{(2)}-\lambda I$ is nonsingular.
 \end{lem}
\begin{pf} Without loss of generality we let $s_2=1$ and write
\begin{equation}\label{def:Aa}
\mathcal{L}^{(2)}=\left[\begin{array}{cc}\ell_{11}&\bi{l}_{12}\\\bi{l}_{21} &L_{22}\end{array}\right]
\end{equation}
where $L_{22}$ is a principal submatrix of $\mathcal{L}^{(2)}$ of order $k_2-1$ and $\bi{l}_{12}^T=\bi{l}_{21}$.
The nonzeroness of the first entry of every eigenvector of $\mathcal{L}^{(2)}$ implies the nonorthogonality of $\bi{l}_{21}$ to the eigenvectors of $L_{22}$~\cite{agu15}, meaning that $\mathcal{L}^{(2)}$ and $L_{22}$ have different eigenvalues~\cite{hw04}. In case the eigenvalue $\lambda$ of $\mathcal{G}_1(\mathbb{G}_2,s_2)$ is also an eigenvalue of $\mathcal{L}^{(2)}$, then $\lambda$ is not an eigenvalue of $L_{22}$ and thus $L_{22}-\lambda I$ are nonsingular. If $\lambda$ is not an eigenvalue of $\mathcal{L}^{(2)}$, $\mathcal{L}^{(2)}-\lambda I$ is nonsingular, meaning that all rows of
$\mathcal{L}^{(2)}-\lambda I$ are independent. If we replace its first row with $\bi{e}_1^T=[\,1\,0\,0\,\cdots\,0\,]$ to yield
\begin{equation}\label{def:Aa1}
\hat{\mathcal{L}}^{(2)}(\lambda):=\left[\begin{array}{cc}1&\\\bi{l}_{21} &L_{22}-\lambda I\end{array}\right],
\end{equation}
then $\hat{\mathcal{L}}^{(2)}(\lambda)$ is nonsingular. To see this, consider $\mathcal{L}^{(1,2)}_{s_2}-\lambda I$. If $\hat{\mathcal{L}}^{(2)}(\lambda)$ is singular, its first row must be a linear combination of other rows, and thus the rank of $\mathcal{L}^{(1,2)}_{s_2}-\lambda I$ is the rank of $I_{k_1}\otimes(\mathcal{L}^{(2)}-\lambda I)$, which is $k_1\times k_2$, a contradiction since $\lambda$ is an eigenvalue of $\mathcal{L}^{(1,2)}_{s_2}$. Note that $|\hat{\mathcal{L}}^{(2)}(\lambda)|=|L_{22}-\lambda I|$, we conclude that $L_{22}-\lambda I$ is nonsingular.
\end{pf}
\begin{thm}\label{thm:idctr}
Suppose $p$ is some positive integer and in every eigenvector of $\mathcal{L}^{(1)}$, there are $p$ fixed positions in which the entries are nonzero. If the $s_2$th entry of every eigenvector of $\mathcal{L}^{(2)}$ is nonzero, then the eigenvalue $\lambda$ of $\mathcal{L}^{(1,2)}_{s_2}$ is not repeated. In particular, for every eigenvector of $\mathcal{L}^{(1,2)}_{s_2}$, the $p$ entries, say the $i_1$th, $i_2$th,$\cdots$, and the $i_p $th entries , corresponding to the $p$ composite vertices of $\mathcal{G}_1(\mathbb{G}_2,s_2)$ are nonzero.
\end{thm}
\begin{pf} We assume $s_2=1$ to simplify the presentation. Instead of considering $\bi{v}:=[\,v_1\,\cdots,v_{k_1k_2}\,]$ in the null space of $\mathcal{L}^{(1,2)}_{s_2}-\lambda I$, we permutate the entries of $\mathcal{L}^{(1,2)}_{s_2}$ and $\bi{v}$ to study $\tilde{\mathcal{L}}^{(1,2)}_{\tilde{s}_2}$ and $\tilde{\bi{v}}:=\left[\,\bar{\bi{v}}\,\hat{\bi{v}}\,\right]$ where
\begin{equation}
\begin{split}\label{def:vbh}
\bar{\bi{v}}&:=\left[\,v_1\,v_{k_2+1}\,\cdots\,v_{(k_1-1)k_2+1}\,\right],\\
\hat{\bi{v}}&:=\left[\,v_2\,v_3\,\cdots\,v_{k_2}\,v_{k_2+2}\,v_{k_2+3}\,\cdots\,v_{2k_2}\,v_{2k_2+2}\cdots\,v_{k_1k_2}\right].
\end{split}
\end{equation}
The nonsingularity of $L_{22}-\lambda I$ by Lemma~\ref{lem:nsing} implies that $\bi{v}$ is in the null space of $\mathcal{L}^{(1,2)}_{s_2}-\lambda I$ if and only if $\tilde{\bi{v}}$ is in the null space of the matrix
\begin{align*}
\tilde{\mathscr{L}}(\lambda):&=\tilde{\mathcal{L}}^{(1,2)}_{\tilde{s}_2}-\lambda I\\
&=\left[\begin{array}{ll}\mathscr{L}_{11}(\lambda) &\\
\mathscr{L}_{21} &\mathscr{L}_{22}(\lambda)\end{array}\right]
\end{align*}
where
\begin{equation}
\begin{split}
\mathscr{L}_{11}(\lambda):&=\mathcal{L}^{(1)}-\left(\lambda+\bi{l}_{21}^T(L_{22}-\lambda I)^{-1}\bi{l}_{21}\right)I,\\
\mathscr{L}_{22}(\lambda):&=I_{k_1}\otimes(L_{22}-\lambda I),\\
\mathscr{L}_{21}:&=I_{k_1}\otimes\bi{l}_{21}.
\end{split}
\end{equation}
This implies that $\bar{\bi{v}}$ is in the null space of $\mathscr{L}_{11}(\lambda)$.
Note that $\tilde{\mathscr{L}}(\lambda)$ is block-wise lower triangular and $\mathscr{L}_{22}(\lambda)$ is nonsingular by Lemma~\ref{lem:nsing}. If $\lambda$ is repeated, the rank of $\tilde{\mathcal{L}}(\lambda)$ is at most $k_1+k_2-2$ and thus the rank of $\mathscr{L}_{11}(\lambda)$ is at most $k_1-2$. This is contradictory to the fact by Lemma~\ref{lem:indc} that the rank of $\mathscr{L}_{11}(\lambda)$ is at least $k_1-1$. We conclude that $\mathcal{L}^{(1,2)}_{s_2}$ does not have a repeated eigenvalue. It also follows from Lemma~\ref{lem:indc} that for all eigenvectors of $\mathcal{L}^{(1,2)}_{s_2}$, the $p$ entries corresponding to their composite vertices of $\mathcal{G}_1(\mathbb{G}_2,s_2)$ cannot be zero.
\end{pf}

\begin{thm}\label{thm:ctr} Following Theorem~\ref{thm:idctr}, $\mathbb{G}_1$ is Laplacian controllable by an input connected to one of its vertices if and only if $\mathcal{G}_1(\mathbb{G}_2,s_2)$ is Laplacian controllable by the input connected to that vertex.
\end{thm}
\begin{pf}
The set of eigenvalues of $\mathcal{L}_{s_2}^{(1,2)}$ can be derived from the set of eigenvalues of $\mathcal{L}^{(2)}+\lambda_i\bi{e}_{s_2}\bi{e}_{s_2}^T$ where $\lambda_i$ is the $i$th
smallest eigenvalue of $\mathcal{L}^{(1)}$, $i\in\{1,2,\cdots,k_1\}$. Suppose $(\lambda_i,\bi{v}^{(i)})$ is an eigenpair of $\mathcal{L}^{(1)}$, and $U^{(i)}$ is an orthogonal modal matrix that diagonalizes $\mathcal{L}^{(2)}+\lambda_i\bi{e}_{s_2}\bi{e}_{s_2}^T$. Then, an orthogonal set of eigenvectors of $\mathcal{L}_{s_2}^{(1,2)}$ is
\begin{equation}\label{eq:evs}
\left\{\,\bi{v}^{(1)}\otimes U^{(1)},
\bi{v}^{(2)}\otimes U^{(2)},\cdots,\bi{v}^{(k_1)}\otimes U^{(k_1)}\,\right\}
\end{equation}
In light of Theorem~\ref{thm:idctr}, it remains to show the case that $\mathbb{G}_1$ is not Laplacian controllable by the input connected to that particular vertex. In this case the (composite) vertex of $\mathbb{G}_1$ induces a zero entry in some $\bi{v}^{(i)}$ where $i\in\{1,2,\cdots,k_1\}$. Thus the same vertex induces a zero entry in some eigenvector in~(\ref{eq:evs}), which leads to the Laplacian uncontrollability of the composite graph by the input.
\end{pf}

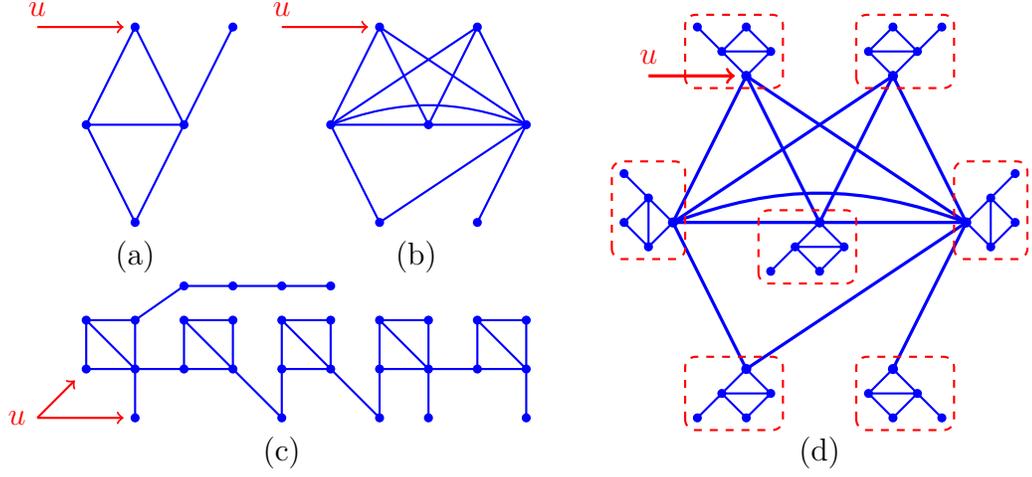
\begin{figure*}[t!]
		\centering
		\begin{tikzpicture}[blue, thick, scale=0.65]
		\newcommand*{\lx}{1}
		\newcommand*{\ly}{3}
		\newcommand*{\lz}{13}
		\newcommand*{\lw}{6}
		\newcommand*{\lv}{0.5}
		\foreach \v in {0,2,4,6,8}{
			\filldraw { (\lx,-2*\lx) ++ (\v*\lx,0) } circle (2pt) -- ++(0:\lx) circle (2pt);
			\filldraw { (\lx,-3*\lx) ++ (\v*\lx,0) } circle (2pt) -- ++(0:\lx) circle (2pt);
			\draw { (\lx,-3*\lx) ++ (\v*\lx,0) } -- ++(90:\lx) -- ++(\lx,-\lx) -- ++(90:\lx);
		}
		\filldraw (2*\lx,-2*\lx) -- ++(\lx,0.7*\lx) circle (2pt) -- ++(\lx,0) circle (2pt) -- ++(\lx,0) circle (2pt) -- ++(\lx,0) circle (2pt);
		\filldraw (2*\lx,-4*\lx) circle (2pt) -- ++(90:\lx);
		\draw (2*\lx,-3*\lx) -- ++(0:\lx);
		\filldraw (5*\lx,-4*\lx) circle (2pt);
		\filldraw (7*\lx,-4*\lx) circle (2pt);
		\filldraw (8*\lx,-4*\lx) circle (2pt) -- ++(90:\lx);
		\filldraw (10*\lx,-4*\lx) circle (2pt) -- ++(90:\lx);
		\draw (4*\lx,-3*\lx) -- ++(\lx,-\lx) -- ++(90:\lx);
		\draw (6*\lx,-3*\lx) -- ++(\lx,-\lx) -- ++(90:\lx);
		\draw (8*\lx,-3*\lx) -- ++(\lx,0);
		\node (1a) at (\lx,-3*\lx) {};
		\node (1b) at (2*\lx,-4*\lx) {};
		\draw [red, ->] (0,-4*\lx) node [left] {$u$} -- (1a);
		\draw [red, ->] (0,-4*\lx) -- (1b);
		\path (5*\lx,-4.7*\lx) node [black] {(c)};
		\filldraw (\lx,2*\lx) circle (2pt) -- ++(2*\lx,0) circle (2pt) -- ++(\lx,2*\lx) circle (2pt);
		\filldraw (\lx,2*\lx) -- ++(\lx, 2*\lx) circle (2pt) -- ++(\lx,-2*\lx);
		\filldraw (\lx,2*\lx) -- ++(\lx,-2*\lx) circle (2pt) -- ++(\lx, 2*\lx);
		\node (2) at (2*\lx,4*\lx) {};
		\draw [red, ->] (0,4*\lx) node [above] {$u$} -- (2);
		\path (2*\lx,-0.7*\lx) node [black] {(a)};
		\filldraw { (0,2*\lx) ++ (\lw,0) } circle (2pt) -- ++(2*\lx,0) circle (2pt) -- ++(2*\lx,0) circle (2pt) -- ++(-\lx,-2*\lx) circle (2pt);
		\filldraw { (0,2*\lx) ++ (\lw,0) } -- ++(  \lx, 2*\lx) circle (2pt) -- ++(3*\lx,-2*\lx);
		\filldraw { (0,2*\lx) ++ (\lw,0) } -- ++(3*\lx, 2*\lx) circle (2pt) -- ++(  \lx,-2*\lx);
		\filldraw { (0,2*\lx) ++ (\lw,0) } -- ++(  \lx,-2*\lx) circle (2pt) -- ++(3*\lx, 2*\lx);
		\filldraw { (2*\lx,2*\lx) ++ (\lw,0) } -- ++(-\lx,2*\lx);
		\filldraw { (2*\lx,2*\lx) ++ (\lw,0) } -- ++( \lx,2*\lx);
		\draw (\lw,2*\lx) edge [out= 20, in =160] ++(4*\lx,0);
		\node (3) at (\lx+\lw,4*\lx) {};
		\draw [red, ->] (-\lx+\lw,4*\lx) node [above] {$u$} -- (3);
		\path (\lw+1.75*\lx,-0.7*\lx) node [black] {(b)};
		\filldraw [very thick] { (0,0) ++ (\lz,0) } circle (2pt) -- ++(\ly,0) circle (2pt) -- ++(\ly,0) circle (2pt) -- ++(-0.5*\ly,-\ly) circle (2pt);
		\filldraw [very thick] { (0,0) ++ (\lz,0) } -- ++ (0.5*\ly, \ly) circle (2pt) -- ++(1.5*\ly,-\ly);
		\filldraw [very thick] { (0,0) ++ (\lz,0) } -- ++ (1.5*\ly, \ly) circle (2pt) -- ++(0.5*\ly,-\ly);
		\filldraw [very thick] { (0,0) ++ (\lz,0) } -- ++ (0.5*\ly,-\ly) circle (2pt) -- ++(1.5*\ly, \ly);
		\filldraw [very thick] { (\ly,0) ++ (\lz,0) } -- ++(-0.5*\ly,\ly);
		\filldraw [very thick] { (\ly,0) ++ (\lz,0) } -- ++( 0.5*\ly,\ly);
		\draw [very thick] (\lz,0) edge [out= 20, in =160] ++(2*\ly,0);
		\filldraw { (0,0) ++ (\lz,0) } -- ++(-\lv, \lv) circle (2pt) -- ++(-\lv,\lv) circle (2pt);
		\filldraw { (0,0) ++ (\lz,0) } -- ++(-\lv,-\lv) circle (2pt) -- ++(-\lv,\lv) circle (2pt);
		\draw { (-2*\lv,0) ++ (\lz,0) } -- ++( \lv, \lv) -- ++(270:2*\lv);
		\draw [red, rounded corners, dashed] { (-2.5*\lv,-1.5*\lv) ++ (\lz,0) } rectangle ++(3*\lv,4*\lv);
		\filldraw { (\ly,0) ++ (\lz,0) } -- ++(-\lv,-\lv) circle (2pt) -- ++(-\lv,-\lv) circle (2pt);
		\filldraw { (\ly,0) ++ (\lz,0) } -- ++( \lv,-\lv) circle (2pt) -- ++(-\lv,-\lv) circle (2pt);
		\draw { (\ly,-2*\lv) ++ (\lz,0) } -- ++(-\lv, \lv) -- ++(2*\lv,0);
		\draw [red, rounded corners, dashed] { (\ly-2.5*\lv,-2.5*\lv) ++ (\lz,0) } rectangle ++(4*\lv,3*\lv);
		\filldraw { (2*\ly,0) ++ (\lz,0) } -- ++( \lv, \lv) circle (2pt) -- ++( \lv, \lv) circle (2pt);
		\filldraw { (2*\ly,0) ++ (\lz,0) } -- ++( \lv,-\lv) circle (2pt) -- ++( \lv, \lv) circle (2pt);
		\draw { (2*\ly+2*\lv,0) ++ (\lz,0) } -- ++(-\lv, \lv) -- ++(270:2*\lv);
		\draw [red, rounded corners, dashed] { (2*\ly-0.5*\lv,-1.5*\lv) ++ (\lz,0) } rectangle ++(3*\lv,4*\lv);
		\filldraw { (0.5*\ly,\ly) ++ (\lz,0) } -- ++(-\lv, \lv) circle (2pt) -- ++(-\lv, \lv) circle (2pt);
		\filldraw { (0.5*\ly,\ly) ++ (\lz,0) } -- ++( \lv, \lv) circle (2pt) -- ++(-\lv, \lv) circle (2pt);
		\draw { (0.5*\ly,\ly+2*\lv) ++ (\lz,0) } -- ++(-\lv,-\lv) -- ++(2*\lv,0);
		\draw [red, rounded corners, dashed] { (0.5*\ly-2.5*\lv,\ly-0.5*\lv) ++ (\lz,0) } rectangle ++(4*\lv,3*\lv);
		\filldraw { (1.5*\ly,\ly) ++ (\lz,0) } -- ++( \lv, \lv) circle (2pt) -- ++( \lv, \lv) circle (2pt);
		\filldraw { (1.5*\ly,\ly) ++ (\lz,0) } -- ++(-\lv, \lv) circle (2pt) -- ++( \lv, \lv) circle (2pt);
		\draw { (1.5*\ly,\ly+2*\lv) ++ (\lz,0) } -- ++( \lv,-\lv) -- ++(180:2*\lv);
		\draw [red, rounded corners, dashed] { (1.5*\ly-1.5*\lv,\ly-0.5*\lv) ++ (\lz,0) } rectangle ++(4*\lv,3*\lv);
		\filldraw { (0.5*\ly,-\ly) ++ (\lz,0) } -- ++(-\lv,-\lv) circle (2pt) -- ++(-\lv,-\lv) circle (2pt);
		\filldraw { (0.5*\ly,-\ly) ++ (\lz,0) } -- ++( \lv,-\lv) circle (2pt) -- ++(-\lv,-\lv) circle (2pt);
		\draw { (0.5*\ly,-\ly-2*\lv) ++ (\lz,0) } -- ++(-\lv, \lv) -- ++(2*\lv,0);
		\draw [red, rounded corners, dashed] { (0.5*\ly-2.5*\lv,-\ly-2.5*\lv) ++ (\lz,0) } rectangle ++(4*\lv,3*\lv);
		\filldraw { (1.5*\ly,-\ly) ++ (\lz,0) } -- ++( \lv,-\lv) circle (2pt) -- ++( \lv,-\lv) circle (2pt);
		\filldraw { (1.5*\ly,-\ly) ++ (\lz,0) } -- ++(-\lv,-\lv) circle (2pt) -- ++( \lv,-\lv) circle (2pt);
		\draw { (1.5*\ly,-\ly-2*\lv) ++ (\lz,0) } -- ++( \lv, \lv) -- ++(-2*\lv,0);
		\draw [red, rounded corners, dashed] { (1.5*\ly-1.5*\lv,-\ly-2.5*\lv) ++ (\lz,0) } rectangle ++(4*\lv,3*\lv);
		\node (5) at (0.5*\ly+\lz,\ly) {};
		\draw [red, ->, very thick] (-\lv+\lz,\ly) node [above] {$u$} -- (5);
		\path (\ly+\lz,-4.7*\lx) node [black] {(d)};
		\end{tikzpicture}
\caption{Subigure (d) is a composite graph whose cell graph and structure graph are a $5$-vertex antiregular graph shown in (a) and a $7$-vertex antiregular graph in (b) respectively. Subfigure (c) extends a $5$-vertex path to a graph that interconnects five $5$-vertex antiregular graphs with the method introduced before (\ref{def:BL}), and connects the resulting graph to a $4$-vertex path shown in the upper part of (c). By Theorem~\ref{thm:ctr} the composite graph in (d) is Laplacian controllable by the input $u$ specified in the subfigure. Lemma~\ref{lem:nzo} and~\ref{lem:apath} suggest that the interconnected graph in (c) is also Laplacian controllable by the specified $u$, with or without the connection to the $4$-vertex path. See Examplpes~\ref{exmp:fig_e} and~\ref{exmp:fig_d} as well.}\label{fig_exmp}
\end{figure*}

\begin{exmp}\label{exmp:fig_e} Subfigure (d) in Figure~\ref{fig_exmp} is a composite graph whose cell graph and structure graph are the $5$-vertex and $7$-vertex  antiregular graphs shown in (a) and (b) respectively. The composite graph is Laplacian controllable by the specified input $u$ since its structure graph is Laplacian controllable by the specified input $u$ shown in (b).
\end{exmp}

Let $j$ be a positive integer and the set $\mathcal{C}_j$ be defined as
\begin{equation}\label{def:cj}
\mathcal{C}_j:=\{\,j,j+(2j+1),j+2(2j+1),j+3(2j+1),\cdots\,\}.
\end{equation}
A special case of Theorem~\ref{thm:idctr} follows.
\begin{cor}
Suppose $\mathbb{G}_1$ in Theorem~\ref{thm:idctr} is a path with $k_1$ vertices and an input is connected to one of these vertices. The input divides the vertices of the path into two parts: one with $k_{11}$ vertices and the other with $k_{12}$ where $k_{11}+k_{12}=k_1-1$. If $\mathbb{G}_2$ and $s_2$ satisfy the conditions in Theorem~\ref{thm:idctr}, then $\mathcal{G}_1(\mathbb{G}_2,s_2)$ is single-input Laplacian controllable if and only if $k_{11}$ and $k_{12}$ do not appear in the same $\mathcal{C}_j$ in~(\ref{def:cj}), $j\in\{1,2,\cdots\}$.
\end{cor}

\begin{pf}
It can be proved that $\mathbb{G}_1$ has a Laplacian eigenvector whose entry corresponding to the vertex connected to the input is zero if and if $k_{11}$ and $k_{12}$ are both in $\mathcal{C}_j$ where $j$ is some positive integer~\cite{hsu17a}. This can also be seen from~\cite[Theorem~5]{yueh05} that for a set of orthogonal Laplacian eigenvectors $\bi{v}^{(1)}, \bi{v}^{(2)},\cdots,\bi{v}^{(k_1)}$ of $\mathbb{G}_1$, $v^{(i)}_j$, the $j$th entry of $\bi{v}^{(i)}$, can be written as
\begin{equation}
v^{(i)}_j =\cos\frac{(i-1)(2j-1)\pi}{2k_1}
\end{equation}
where $i,j\in\{1,2,\cdots,k_1\}$. If $k_1$ and $k_2$ are in $\mathcal{C}_{j^*}$, we can write $k_{11}=j^*+n_1(2j^*+1)$, $k_{12}=j^*+n_2(2j^*+1)$, where $n_1,n_2$ are some nonnegative integers, and thus $k_1=(n_1+n_2+1)(2j^*+1)$. The entry corresponding to the vertex connected to the input can be written as
\begin{equation}
v^{(i)}_{k_{11}+1} =\cos\left\{(i-1)\frac{2n_1+1}{n_1+n_2+1}\frac{\pi}{2}\right\},
\end{equation}
for $i\in\{1,2,\cdots,(n_1+n_2+1)(2j^*+1)\}$. Thus
\begin{equation}
v^{(n_1+n_2+2)}_{k_{11}+1} =\cos\left\{(2n_1+1)\frac{\pi}{2}\right\}=0.
\end{equation}
For the same corresponding vertex, Theorem~\ref{thm:idctr} shows that the nonzero property of some entries in Laplacian eigenvectors of $\mathbb{G}_1$ is preserved in those of $\mathcal{G}_1(\mathbb{G}_2,s_2)$. It can also be seen that if some vertex induces a zero entry in some Laplacian eigenvector of $\mathbb{G}_1$, it will also induce a zero entry in some Laplacian eigenvector of $\mathcal{G}_1(\mathbb{G}_2,s_2)$. Thus the proof is completed.
\end{pf}

\subsection{Interconnecting antiregular graphs}
We have proposed in the previous section a condition for a class of single-input Laplacian controllable graph to maintain its single-input controllability after interconnecting a finite number of its identical copies. This condition requires the interconnection via the composite vertices. In the following we show that this requirement can be relaxed by interpreting the connection structure of a path from a novel perspective. That is, a path can be regarded as a graph generated by interconnecting several two-vertex antiregular graphs and possibly connected to a one-vertex path. Motivated by this interpretation, we propose the second interconnection scheme to preserve the single-input Laplacian controllability.

Let $S$ be a real and symmetric matrix of order $k_2$ and
 \begin{equation}\label{def:ms}
\mathcal{S}:=I_{k_1}\otimes S+\sum_{i=1}^{k_1-1}\bi{z}_i\bi{z}_i^T
\end{equation}
 where $\bi{z}_i:=c_1\bi{e}_{(i-1)k_2+p}+c_2\bi{e}_{ik_2+1}$, $p\in\{1,2,\cdots,k_2\}$, and $c_1,c_2$ nonzero.

\begin{figure*}[t!]
{\footnotesize
\begin{equation}\label{def:mathA1}
\tilde{\mathcal{S}}(\lambda):=\left[\begin{array}{cccccccccc}
0&\cdots &0 &\cdots &0&&\cdots&\cdots&0&\cdots\\
\bi{s}_{21} & S_{22}-\lambda I &&&&&&&\\
c_1c_2 & &c_2^2 &&&&&&\\
& &\bi{s}_{21} &S_{22}-\lambda I &&&&&\\
& & c_1c_2 & &c_2^2 &&&&&\\
&&&\ddots& &&\ddots &&&\\
&&&&&&\bi{s}_{21}&S_{22}-\lambda I&&\\
&&&&&&c_1c_2&&c_2^2&\\
&&&&&&&&\bi{s}_{21} &S_{22}-\lambda I \\
\end{array}\right].
\end{equation}
}
\end{figure*}

\begin{lem} \label{lem:vuni}
If in~(\ref{def:ms}) the first entry of every eigenvector of $S$ is nonzero, then every eigenvalue $\lambda$ of $S$ is a distinct eigenvalue of $\mathcal{S}$ and the first entry of $\bi{v}$ is nonzero where $(\lambda,\bi{v})$ is an eigenpair of $\mathcal{S}$. If, in addition, the $p$th entry of every eigenvector of $S$ is also nonzero, then so are the $((i-1)k_2+1)$th entries, for every $i\in\{2,3,\cdots,k_1\}$, of the corresponding eigenvectors of $\mathcal{S}$.
\end{lem}

\begin{pf}
Write
\begin{equation}
S:=\left[\begin{array}{cc}s_{11}&\bi{s}_{12}\\\bi{s}_{21} &S_{22}\end{array}\right]
\end{equation}
where $S_{22}$ is a principal submatrix of $S$ of order $k_2-1$ and $\bi{s}_{12}^T=\bi{s}_{21}$. Since the first entry of every eigenvector of $S$ is nonzero, the eigenvalues of $S$ are distinct~\cite{rah09}. Let $(\lambda_i,\bi{v}_i)$ be an eigenpair of $S$ for every $i\in\{1,2,\cdots,k_2\}$ where $\bi{v}_i=[\,v_{i1}\,v_{i2}\,\cdots\,v_{ik_2}\,]^T$ and $v_{i1}\neq 0$. An apparent eigenpair of $\mathcal{S}$ is
$(\lambda_i,\bi{u}_i)$ where $\bi{u}_i:=[\,1\,t\,t^2\cdots\,t^{k_1-1}]^T\otimes\bi{u}_i^T$ with
$t=-(c_1v_{ip})/(c_2v_{i1})$ due to the orthogonality of $\bi{z}_i$ to $\bi{u}_j$ for every $i\in\{1,2,\cdots,k_1-1\}$ and $j\in\{1,2,\cdots,k_2\}$. Clearly the first entry of $\bi{u}_i$ is nonzero. If $v_{ip}\neq 0$ then
the $((j-1)k_2+1)$th entry, for every $j\in\{2,3,\cdots,k_1\}$ of $\bi{u}_i$, is also nonzero. Now we show that $\lambda$ is a unique eigenvalue of $\mathcal{S}$ if it is an eigenvalue of $S$. We first consider the case that $p=1$. The result in this case follows from Theorem~\ref{thm:idctr}. However we provide a different proof to make the argument applicable to the cases not covered by Theorem~\ref{thm:idctr}. Similar to the reason in proving Lemma~\ref{lem:nsing}, the nonzeroness of the first entry of every eigenvector of $S$ suggests that $S$ and $S_{22}$ do not share a common eigenvalue. After some row operations, we can write the null space of $\mathcal{S}-\lambda I$ as that of $\tilde{\mathcal{S}}(\lambda)$ defined in~(\ref{def:mathA1}), whose first row is a zero row. If $\bi{u}^T\tilde{\mathcal{S}}(\lambda)=\textbf{0}^T$ where $\bi{u}=[\,u_1\,u_2\,\cdots\,u_{k_1k_2}\,]$, the row independence of $S_{22}-\lambda I$ implies that $u_{(k_1-1)k_2+i}=0$ for every $i\in\{2,3,\cdots,k_2\}$, and thus $u_{(k_1-1)k_2+1}=0$. In fact, from the structure of $\tilde{\mathcal{S}}(\lambda)$ we obtain that $u_j=0$ for $j\in\{2,3,\cdots,k_1k_2\}$. That is, the rank of $\tilde{\mathcal{S}}(\lambda)$ is $k_1k_2-1$ as $\lambda$ is an eigenvalue of $S$. If $p\neq 1$, the only differences in $\tilde{\mathcal{S}}(\lambda)$ are the positions of $c_1c_2$. Similar arguments can be applied to prove that the rank of $\tilde{\mathcal{S}}(\lambda)$ is $k_1k_2-1$ and thus $\lambda$ is not repeated.
\end{pf}

Now we consider a different interconnection scheme from that used to generate a composite graph. Let
$\mathbb{G}^{(k)}_A$ denote a $k$-vertex antiregular graph and its Laplacian matrix $\mathcal{L}^{(k)}_A$, namely,
\begin{equation*}
-\begin{bmatrix}
-(k-1) &    1   & \cdots &    1   &    1   & \cdots &    1   &    1   \\
1   & -(k-2) & \cdots &    1   &    1   & \cdots &    1   &        \\
\vdots & \vdots & \ddots & \vdots & \vdots &\iddots        &        &        \\
1   &    1   & \cdots & -\lfloor k/2\rfloor & \beta_k & & & \\
1   &    1   & \cdots & \beta_k  & -\lfloor k/2\rfloor & & & \\
\vdots & \vdots &\iddots        &        &        & \ddots &        &        \\
1   &    1   &        &        &        &        &   -2   &        \\
1   &        &        &        &        &        &        &   -1
\end{bmatrix}
\end{equation*}
where $\beta_k$ is $1$ if $k$ is even and is $0$ otherwise.
Note that an antiregular graph is also a threshold graph. Its Laplacian spectrum can be readily derived from its degree conjugate. That is, the set of eigenvalues of $\mathcal{L}_A^{(k)}$ is
$\{0,1,\cdots,k\}\setminus\{\lceil k/2\rceil\}$.
A full set of orthogonal eigenvectors of $\mathcal{L}^{(k)}_A$ can be obtained using the following lemma.
\begin{lem}\cite{hsu16,agu15b} \label{lem:ev}
	Let the $(i,j)$th entry of matrices $T^{(m)}$ be $t_{ij}^{(m)}$ for each $m\in\{1,2,3,4\}$. Suppose $T^{(1)}=\mathcal{L}^{(k)}_A$ and let $T^{(2)},T^{(3)}$ be generated by
	\begin{equation} \label{def:t2}
	t_{ij}^{(2)}=\left\{
	\begin{array}{cl}
	-1-t_{ij}^{(1)}, & \mbox{if $ j>i$} \\
	t_{ij}^{(1)},    & \mbox{o.w.}
	\end{array}
	\right.
	\end{equation}
	and
	\begin{equation} \label{def:t3}
	t_{ij}^{(3)}=\left\{
	\begin{array}{cl}
	-\sum_{k,k\ne j} t_{kj}^{(2)}, & \mbox{if  $j=i$} \\
	t_{ij}^{(2)},                  & \mbox{o.w.}
	\end{array}
	\right..
	\end{equation}
Finally, remove the (unique) zero column of $T^{(3)}$ and append the column of $1$'s (or $-1$'s) to the last column to yield $T^{(4)}$. Then the $j$th column of $T^{(4)}$ is the eigenvector corresponding to the $j$th largest eigenvalue of $\mathcal{L}^{(k)}_A$, or the $j$th entry of the conjugate of the degree sequence of $\mathbb{G}^{(k)}_A$.
\end{lem}

To interconnect a finite number of identical antiregular graphs, we start with one $k_2$-vertex antiregular graph, and repeatedly add in a new $k_2$-vertex antiregular graph. When the $m$th antiregular graph is added, $m=2,3,\cdots$, one of its two degree-repeating vertices is connected to the terminal or dominating vertex of the $(m-1)$th antiregular graph that was added. Let $\mathbb{G}^{(1,2)}_A$ be the resulting graph of interconnecting $k_1$ antiregular graphs, each with $k_2$ vertices. The corresponding Laplacian matrices $\mathcal{L}_A^{(1,2)}$ can be written as
\begin{equation}\label{def:BL}
\mathcal{L}_A^{(1,2)}:=I_{k_1}\otimes\mathcal{L}_A^{(k_2)}+\sum_{i=1}^{k_1-1}\bi{z}_i\bi{z}_i^T,\\
\end{equation}
where$\bi{z}_i\in\left\{\bi{e}_{(i-1)k_2+1}-\bi{e}_{ik_2+\bar{\kappa}_2}, \;\bi{e}_{ik_2}-\bi{e}_{ik_2+\bar{\kappa}_2}\right\}$
and $\bar{\kappa}_i=\lceil k_i/2\rceil$, $i\in\{1,2\}$. For integers $a,b$ with $a<b$, let
\begin{equation}
\mathcal{R}^a_b:=R^{(-1)}_{b,b-1}\cdots R^{(-1)}_{a+2,a+1}R^{(-1)}_{a+1,a}
\end{equation}
where $R_{\alpha,\beta}^{(\gamma)}$ is the elementary row operation matrix, namely, the identity matrix except its $(\beta,\alpha)$th entry being $\gamma$. As a result,
\begin{align*}
\mathscr{L}_A^{(1,2)}(\lambda):&=\mathcal{R}_{k_1k_2}^{(k_1-1)k_2+1}\cdots
\mathcal{R}_{2k_2}^{k_2+1}\mathcal{R}_{k_2}^1\left(\mathcal{L}_A^{(1,2)}-\lambda I\right)\\
&=I_{k_1}\otimes\tilde{\mathcal{L}}_A^{(k_2)}(\lambda)+\prod_{i=0}^{k_1-1}\mathcal{R}_{(k_1-i)k_2}^{(k_1-i-1)k_2+1}
\sum_{i=1}^{k_1-1}\bi{z}_i\bi{z}_i^T\\
&=I_{k_1}\otimes\tilde{\mathcal{L}}_A^{(k_2)}(\lambda)+\sum_{i=1}^{k_1-1}\tilde{\bi{z}}_i\bi{z}_i^T
\end{align*}
where $\tilde{\mathcal{L}}_A^{(k_2)}$ is shown in~(\ref{eq:ltilde}) and
\begin{equation}
\tilde{\bi{z}}_i=\left\{\begin{array}{ll}\bi{z}_i+\bi{e}_{ik_2+\bar{\kappa}_2-1}&\mbox{if $\bi{z}_i=\bi{e}_{(i-1)k_2+1}-\bi{e}_{ik_2+\bar{\kappa}_2}$},\\
\bi{z}_i-\bi{e}_{ik_2-1}+\bi{e}_{ik_2+\bar{\kappa}_2-1}&\mbox{if $\bi{z}_i=\bi{e}_{ik_2}-\bi{e}_{ik_2+\bar{\kappa}_2}$}.\\
\end{array}\right.
\end{equation}
\begin{figure*}[t!]
\begin{equation}\label{eq:ltilde}
{\footnotesize
\tilde{\mathcal{L}}_A^{(k_2)}(\lambda)=-\left[\begin{array}{rrrrrrrrr}
-k_2+\lambda	&k_2-1-\lambda	 &	&	&		&	&	&	&1\\
 	   &-k_2+1+\lambda &k_2-2-\lambda	 	& 	& 	& 	& 	&1	&  \\
 	   & 	&-k_2+2+\lambda	 &\ddots	& 	 	& 	&1	 & 	&  \\
 	 	& 	&  &\ddots &\ddots	&\iddots 	& 	& 	& \\
 	 	& 	& 	& 	&\ddots	&\ddots 	& 	& 	&\\
 	 	& 	& 	&\iddots	&	 &\ddots	&3-\lambda 	& 	&\\
 	 	& 	&1	& 	& 	& 	&-3+\lambda	&2-\lambda 	& \\
 	    &1	& 	& 	& 	& 	& 	&-2+\lambda	&1-\lambda \\
      1& 	& 	& 	& 	& 	& 	&  &-1+\lambda\\
\end{array}
\right].
}
\end{equation}
\end{figure*}
Observe that $\tilde{\mathcal{L}}_A^{(k_2)}(\lambda)$ has a special structure that only two rows have two entries and other rows three entries. More importantly, the positions of nonzero entries in these rows facilitate the identification of nonzero entries of eigenvectors of $\mathcal{L}_A^{(1,2)}$, as shown in the following Lemma.

\begin{lem}\label{lem:nzo}
The $\bar{\kappa}_2$th and $\big( \bar{\kappa}_2 + 1 \big)$th entries of eigenvectors of $\mathcal{L}_A^{(1,2)}$ are nonzero.
\end{lem}
\begin{pf}
By Lemma~\ref{lem:vuni}, every eigenvalue of $\mathcal{L}_A^{(k_2)}$ is a distinct eigenvalue of $\mathcal{L}_A^{(1,2)}$, and the eigenvectors of $\mathcal{L}_A^{(1,2)}$ corresponding to these eigenvalues have their $\bar{\kappa}_2$th and $\big( \bar{\kappa}_2 + 1 \big)$th entries nonzero. Now we consider the eigenvectors corresponding to the $\lambda$'s that are not the eigenvalues of $\mathcal{L}_A^{(k_2)}$.  The $\bar{\kappa}_2$th row of $\tilde{\mathcal{L}}_A^{(k_2)}(\lambda)$ in~(\ref{eq:ltilde}) implies that the $\bar{\kappa}_2$th and $\big( \bar{\kappa}_2 + 1 \big)$th entries of the eigenvectors have the same value, say $c$. The positions of nonzero entries in the rows of $\tilde{\mathcal{L}}_A^{(k_2)}(\lambda)$ imply that the first $k_2$ entries of the eigenvectors can be written as $[\,(1-\lambda)c\;c\;c\,\cdots\,c\,]$. If $c$ is zero, the first row of $\mathscr{L}_A^{(1,2)}$ suggests that the $(k_2+\bar{\kappa}_2 )$th entry of the eigenvector is also zero. The special structure of $\tilde{\mathcal{L}}_A^{(k_2)}(\lambda)$ implies that the first $2k_2$ entries of the eigenvectors are zero. Continuing this argument yields that the eigenvectors are zero vectors, a contradiction. We thus conclude that $c$ is not zero and the proof is completed.
\end{pf}

The nonzero property in Lemma~\ref{lem:nzo} not only implies the distinctness of Laplacian eigenvalues of $\mathbb{G}^{(1,2)}_A$ but also gives a hint on the vertex selection that ensures the single-input Laplacian controllability of the graph. It was shown in~\cite[Theorem~2.1]{far14} that a necessary and sufficient condition for a control vector to render the graph Laplacian controllable can be derived by analyzing eigenvectors of the bordered matrix~\cite[p.26]{horn13} composed of the Laplacian matrix of the graph and the control vector. However, analyzing this bordered matrix is challenging since some nice properties of Laplacian matrices no longer exist. In our case we  relate the selection of control vecetors for $\mathbb{G}^{(1,2)}_A$ to that of $\mathbb{G}^{(k_2)}_A$. We first differentiate the following two cases of $\mathcal{L}_A^{(1,2)}$ in~(\ref{def:BL}).
\begin{equation}
\mathcal{L}_A^{(1,2)}=\left\{\begin{array}{cc}\bar{\mathcal{L}}_A^{(1,2)} &\mbox{if $\bi{z}_1=\bi{e}_1-\bi{e}_{k_2+\bar{\kappa}_2}$ }\\ \undb{\mathcal{L}}_A^{(1,2)}&\mbox{if $\bi{z}_1=\bi{e}_{k_2}-\bi{e}_{k_2+\bar{\kappa}_2}$}\end{array}\right..
\end{equation}

\begin{thm} Suppose $\tilde{\bi{b}}$ is a binary vector of size $k_2$ and $\bar{\bi{b}}:=\left[\,\tilde{\bi{b}}^T\;\bi{0}^T\,\right]^T$ is a zero-padding expansion of $\tilde{\bi{b}}$, with size $k_1k_2$. The following statements are equivalent:
\begin{enumerate}
\item the sum of the $\bar{\kappa}_2$th and $(\bar{\kappa}_2+1)$th entries of $\bi{b}_1$ is 1;
\item $\left(\mathcal{L}^{(k_2)}_A,\tilde{\bi{b}}\right)$ is Laplacian controllable;
\item $\left(\BL,\bar{\bi{b}}\right)$ is Laplacian controllable.
\end{enumerate}
Similarly, let $\undt{\bi{b}}$  be a binary vector of size $k_2-1$ and $\undb{\bi{b}}:=\left[\,\undt{\bi{b}}^T\;\bi{0}^T\,\right]^T$ is a zero-padding expansion of $\undt{\bi{b}}$,
with size $k_1k_2$. The following statements are equivalent:
\begin{enumerate}
\item the sum of the $\bar{\kappa}_2$th and $(\bar{\kappa}_2+1)$th entries of $\undt{\bi{b}}$ is 1;
\item $\left(\mathcal{L}^{(k_2)}_A,\left[\,\undt{\bi{b}}^T\;0\,\right]^T\right)$ is Laplacian controllable;
\item $\left(\UBL,\undb{\bi{b}},\right)$ is Laplacian controllable.
\end{enumerate}
\end{thm}
\begin{pf} Consider the case of $\BL$. The equivalence of the first two conditions is a well-known result that follows directly from Theorem~\ref{thm:ns} and Lemma~\ref{lem:ev}. In proving Theorem~\ref{lem:nzo} we have shown that the first $k_2$ entries of eigenvector of $\BL$ is either an eigenvector of $\mathcal{L}_A^{(k_2)}$, or in
the form of $[\,(1-\lambda)c\;c\;c\cdots\;c\,]$ where $c$ is nonzero and $\lambda$ is not an eigenvalue of $\mathcal{L}_A^{(k_2)}$. This established the equivalence of the three conditions. The difference between the cases of and $\BL$ and $\UBL$ is that if the first $k_2$ entries of eigenvector of $\UBL$ is not an eigenvector of $\mathcal{L}_A^{(k_2)}$, it has the form of $\bi{u}:=[\,(1-\lambda)c\;c\;c\cdots\;c\;f(\lambda)c\,]^T$ where $c$ is nonzero, $f(\lambda)=\lambda^2-k_2\lambda+1$ and $\lambda$ is not an eigenvalue of $\mathcal{L}_A^{(k_2)}$. The existence of $f(\lambda)$ makes it difficult to determine if $\UBL$ has an eigenvalue leading to a zero sum of some entries in $\bi{u}$. We thus assign a zero to the $k_2$th entry of $\undb{\bi{b}}$ to restrict the control vector and avoid the difficulty. The equivalence of the three conditions then follows again from Theorem~\ref{thm:ns} and Lemma~\ref{lem:ev}.
\end{pf}

\subsection{Appending a path}
In the following lemma we show that a single-input Laplacian controllable graph could preserve its single-input controllability after interconnecting a path.
\begin{lem}\label{lem:apath} Let $A_1,A_2$ and $Z$ be square matrices of order $n_1,n_2$ and $n_1+n_2$ respectively, where
$a_{ij}$ and $z_{ij}$, the $(i,j)$th entries of $A_1$ and $Z$ respectively, satisfy
\begin{equation*}
a_{ij}:\left\{\begin{array}{cc}\neq 0 &\mbox{if $j=i+1\,,\forall i\in\{1,2,\cdots,n_1-1\}$};\\
=0 &\mbox{if $j\ge i+2\,,\forall i\in\{1,2,\cdots,n_1-2\}$};
\end{array}\right.
\end{equation*}
and
\begin{equation*}
z_{ij}:\left\{\begin{array}{cl}\neq 0 &\mbox{if $(i,j)=(n_1,n_1+1)$};\\
=0 &\mbox{if $(i,j)\notin\{(n_1,n_1),(n_1+1,n_1),(n_1+1,n_1+1)\}$}.
\end{array}\right.
\end{equation*}
If the first entry of every eigenvector of $A_2$ is nonzero, then so is that of the matrix:
\begin{equation}
\mathcal{A}:=\left[\begin{array}{cc}A_1&\\&A_2\end{array}\right]+Z.
\end{equation}
 \end{lem}
\begin{pf} If the first entry of some eigenvector of $\mathcal{A}$ is zero, then all but the first column of $\mathcal{A}-\lambda I$ are dependent for some $\lambda$. The special structures of $A_1$ and $Z$ then imply that
 all but the first column of $A_2-\lambda I$ are dependent for some $\lambda$. However, this is a contradiction to Lemma~\ref{lem:indc}.
\end{pf}
\begin{exmp}\label{exmp:fig_d} In subfigure (c) of Figure~\ref{fig_exmp} we present an example that generalizes the Laplacian controllability of a $6$-vertex path. It can be seen that (c) is generated by interconnecting five antiregular graphs, each with five vertices as shown in~(a). Either the terminal vertex or the dominating vertex of each antiregular graph is connected to the degree-repeating vertex of the antiregular graph in its right. In particular, the degree-repeating vertex of the leftmost antiregular graph is connected to a $4$-vertex path. Lemma~\ref{lem:nzo} and~\ref{lem:apath} suggest that the resulting graph is Laplacian controllable by the control input $u$ specified in (c), with or without the connection to the $4$-vertex path.
\end{exmp}

\section{Conclusions and Future Works}
We have studied the Laplacian controllability of a novel class of graphs constructed by interconnecting several single-input Laplacian controllable graphs. The concepts of the composite graph, composite vertex, cell graph and structure graph have been introduced, and used to generate the first type of graphs. Suppose the cell graph is
Laplacian controllable by an input connected to the composite vertex. It has been shown that the structure graph is Laplacian controllable by an input connected to some vertex of the graph, if and only if the composite graph is Laplacian controllable by the input connected to that (composite) vertex. Constructing the second type of graphs is motivated by the observation that a path is actually the resulting graph after interconnecting a finite number of two-vertex antiregular graphs and possibly with one more vertex, which can be viewed as a one-vertex path. We have shown that with a similar connection scheme, the single-input Laplacian controllability of a path is preserved in its generalized version as long as the number of vertices of the antiregular graphs and that of the path are any two nonnegative integers but not zero simultaneously. Our results expand the class of single-input Laplacian controllable graphs significantly, and have potential applications to the design of a Laplacian controllable graph under edge constraints~\cite{hsu17b}. Note that a graph could be single-input Laplacian controllable even if its Laplacian eigenvectors do not have the nonzero property used in our work. Exploring the method of interconnecting such graphs to preserve the single-input Laplacian controllability is interesting. Extending our results to the multi-input case will also be the topic of interest.


\section*{References}
\bibliographystyle{elsarticle-num}
\bibliography{thresh}
\end{document}